\documentclass[12pt]{amsart}

\usepackage[pdftex]{graphicx}
\usepackage{placeins}
\usepackage{geometry}
\usepackage{url}
\usepackage{threeparttable}
\usepackage{booktabs}
\usepackage{appendix}
\usepackage{fancyhdr}
\theoremstyle{remark} 
\newtheorem*{Acknowledgements}{Acknowledgements}

\newcommand{\R}{\ensuremath{\mathbb{R}}}
\newcommand{\cF}{\ensuremath{\mathcal{F}}}
\newcommand{\cM}{\ensuremath{\mathcal{M}}}
\newcommand{\cR}{\ensuremath{\mathcal{R}}}
\DeclareMathOperator{\rank}{rank}

\begin{document}

\title{Nudged Elastic Band in Topological Data Analysis}

\author{Henry Adams}
\address{Department of Mathematics, Duke University, NC 27708}
\email{hadams@math.duke.edu}
\author{Atanas Atanasov}
\address{Department of Mathematics, Harvard University, Cambridge, MA 02138}
\email{nasko@math.harvard.edu}
\author{Gunnar Carlsson}
\address{Department of Mathematics, Stanford University, Stanford, CA 94305}
\email{gunnar@math.stanford.edu}

\begin{abstract}
We use the nudged elastic band method from computational chemistry to analyze high-dimensional data. Our approach is inspired by Morse theory, and as output we produce an increasing sequence of small cell complexes modeling the dense regions of the data. We test the method on data sets arising in social networks and in image processing. Furthermore, we apply the method to identify new topological structure in a data set of optical flow patches.
\end{abstract}

\maketitle

\section{Introduction}

The analysis of large sets of high-dimensional data is a fundamental problem for all branches of science and engineering. Regression analysis can be used effectively when the data has a linear or low degree polynomial structure based on a choice of model for the data. However, it is often the case that the data is genuinely nonlinear and that there isn't an obvious choice for how to model it. The purpose of topological data analysis \cite{TopologyData} is to provide methods which produce simple combinatorial representations of the data.

In this paper we construct a cell complex representation in which the cell structure depends on the density of the points in the data set. We adapt the mathematical formalism of Morse theory, which in its idealized form constructs a cell decomposition of a manifold using sublevel sets of a function on the manifold (called the ``Morse function''), to the setting of point clouds, i.e.\ finite sets of points in Euclidean spaces. The specific features of our approach are as follows.

\begin{itemize}
\item{We use a density estimator as our analogue of the Morse function. One could also use other intrinsic functions on the geometry, such as notions of data depth, to obtain other compact representations.}
\item{In order to sample cells from the analogue of the Morse skeleton of the Morse complex, we adapt the nudged elastic band method (NEB) from computational chemistry \cite{JonssonMillsJacobsen, HenkelmanJonsson, HenkelmanUberuagaJonsson, SheppardTerrellHenkelman} to our situation. NEB has been used to study high-dimensional conformation spaces of complicated molecules, and typically uses an energy function as the analogue of the Morse function.}
\item{We produce an increasing sequence of cell complex models. In accordance with the idea of topological persistence, this increasing sequence gives a more accurate representation of the data than the choice of any single complex.}
\item{We currently construct only the one-dimensional skeleton of the cell complex. Producing higher-order cells will require more difficult mathematics, since the minimization problems involved in the construction of such cells are challenging, and are related to minimal surface problems in geometric analysis.}
\end{itemize}

\begin{figure}[h]
\centerline{\includegraphics[width=0.5\textwidth]{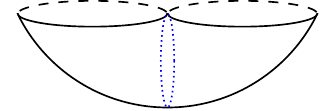}}
\caption{A sublevel set of a torus in $\R^3$ with Morse function given by height. The dotted loop is the Morse skeleton.}
\label{fig: torusSublevel}
\end{figure}

In studying data sets computationally, one finds that outliers will generally obscure topological features. One approach for dealing with outliers is thresholding by density, i.e.\ studying superlevel sets of a function that estimates density. A difficulty is that the superlevel sets of density frequently require large numbers of points to represent them, since they are ``codimension zero" subsets. Consider Figure~\ref{fig: torusSublevel}, which contains a standard Morse theoretic picture: a sublevel set of a torus in $\R^3$ with Morse function given by height. As is standard in Morse theory, the sublevel set is homotopy equivalent to the Morse skeleton, which is the dotted loop in Figure~\ref{fig: torusSublevel}. In order to achieve an accurate representation of the homotopy type, sampling from the entire sublevel set will require many more points than sampling from the Morse skeleton. The goal of this paper is to demonstrate that one can effectively sample cells from the Morse skeleton of a density function on the data set, thereby obtaining a more economical representation of the topology of the data set.

We demonstrate how the method works on several nonlinear data sets. First, while it is common to cluster social network data, there is more understanding to be gained about the structure within a single cluster or about the relationships between different clusters. As illustrated with the social network data set, our method can find such relationships. Second, from a collection of optical image patches we recover a minimal geometric description of the three circle model, which has applications in image compression and in texture analysis \cite{KleinBottle, PereaCarlsson}. Third, we find the primary circle model for range image patches \cite{Range}. Fourth, we identify new topological structure in a data set of optical flow patches.

In Section \ref{S: Related Work} we survey related work, and in Section \ref{S: Background} we provide background material on CW cell complexes, Morse theory, and NEB. We describe our method in Section \ref{S: Method}, present results on data sets in Section \ref{S: Results}, and conclude in Section \ref{S: Conclusion}.

\section{Related Work}\label{S: Related Work}

The analysis of high-dimensional point cloud data is an important yet challenging task for which a wide variety of tools exist. One approach is to map the data points to a lower dimension. This can be done linearly via projection pursuit \cite{Huber} or principal component analysis \cite{Jolliffe}, else nonlinearly via Isomap \cite{TenenbaumDeSilvaLangford}, Locally Linear Embedding \cite{RoweisSaul}, or Laplacian eigenmaps \cite{BelkinNiyogi}, among others. A second approach is to build a combinatorial representation of the data, for example a dendrogram or a cluster tree, and this is the direction we consider.

\begin{figure}[h]
\centerline{\includegraphics[width=0.8\textwidth]{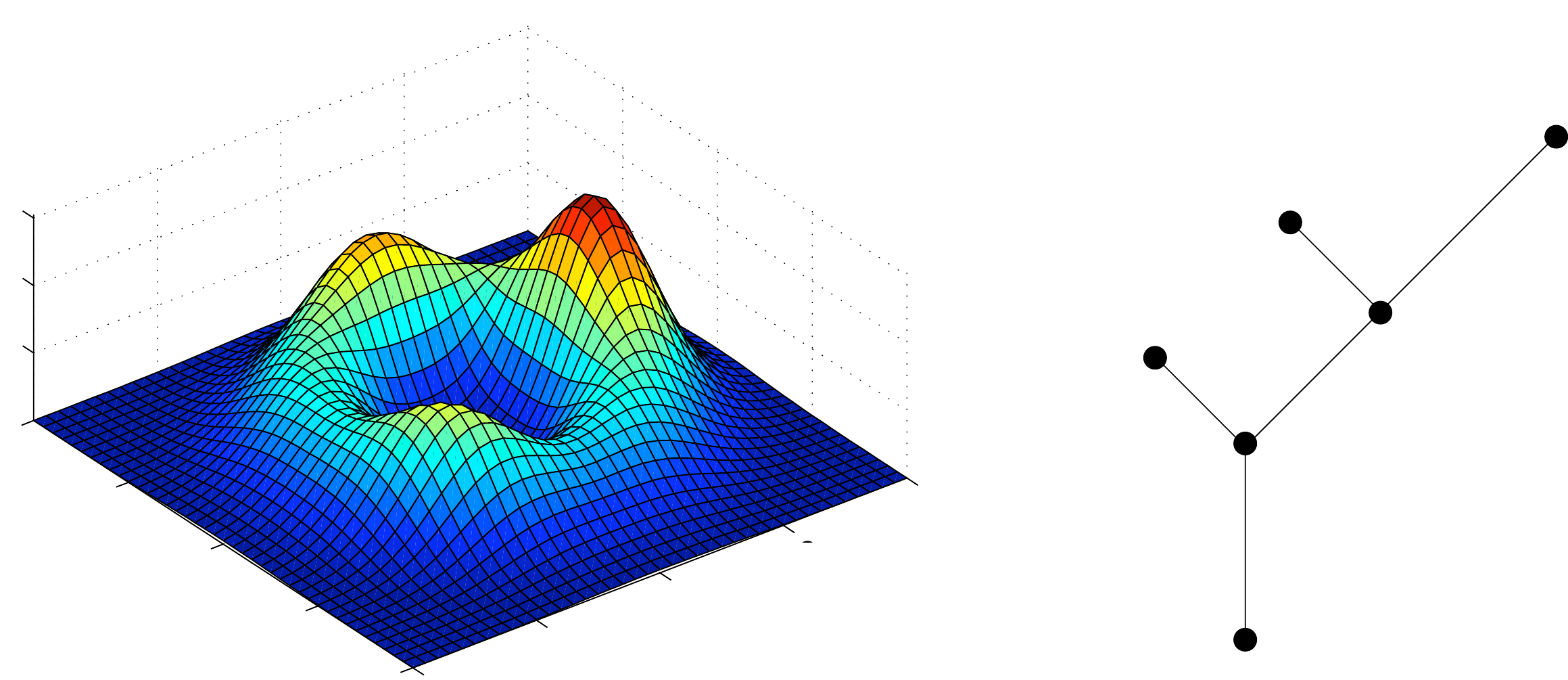}}
\caption{A probability density function on the left and the cluster tree for its superlevel sets on the right. The function has a circular base and three bumps of varying heights. For a certain choice of density threshold the superlevel set of the function has the topology of a circle, but this is not reflected in the cluster tree.}
\label{fig: cluster tree}
\end{figure}

Given a point cloud data set drawn from an unknown probability density function, one may cluster the data by approximating the basins of attraction of the modes of density \cite[pp.\ 205]{Wishart, Hartigan}. The goal of \cite{Stuetzle} is to recover the cluster tree of density, a combinatorial model describing how the connected components of the superlevel sets of density merge as the density threshold decreases. Regardless of its topology, a connected component of a superlevel set is always represented in the cluster tree by a single point. For example, some superlevel sets of the density function in Figure~\ref{fig: cluster tree} have circular topology, but this is not apparent in the cluster tree. There is more information to be gained by considering the topology of the superlevel sets. 

In topological data analysis one often models a superlevel set of density with a simplicial complex, such as a \u Cech, Vietoris--Rips, alpha, or witness complex \cite{AlphaShapes, WitnessComplexes}. These simplicial complexes typically contain thousands of simplices, are too large to interpret by hand, and on their own are not useful descriptors of a superlevel set. Yet from these complexes one can estimate the homology groups of the superlevel set using persistent homology \cite{TopologicalPersistence, ComputingPersistent, ComputationalTopology}. Though homology is a useful signature for describing the geometry of a superlevel set, in general there are several possible models with given homology groups, and determining which model best fits the data is a supervised step. Hence our method, which produces models matching the geometry of the data, plays a complementary role to persistent homology. For example, Carlsson et al.\ in \cite{KleinBottle} use persistent homology to identify an image processing data set with $\rank(H_0) = 1$ and $\rank(H_1) = 5$, where $H_i$ denotes the $i$-th homology group. There are many spaces satisfying these homological constraints, including a wedge sum of five circles. Carlsson et al.\ propose a particular space: a model containing three circles with four intersection points (Figure~\ref{fig: threeCircle}). In Section \ref{SS: Optical image patches} we use our method to obtain this three circle model directly (Figure~\ref{fig: OpticalImagePatches_1cells}). Our model contains only four 0-cells and eight 1-cells, which is much fewer than the number of simplices needed in a \u Cech, Vietoris--Rips, alpha, or witness complex reconstruction. Applications of persistent homology often study only one superlevel set of density at a time due to the difficulty of multidimensional persistent homology \cite{MultidimensionalPersistence}, but exceptions employing kernel density estimators include \cite{ChazalGuibasOudotSkraba, FasyLecciRinaldoWassermanBalakrishnanSingh, BubenikKim}.

There are a variety of adaptations of Morse theory to combinatorial \cite{Forman} or applied settings.
Closely related to the Morse complex that we consider is the Morse-Smale complex, which can be thought of as the intersection of the two Morse complexes for Morse function $f$ and $-f$. The Morse-Smale complex can be efficiently computed for a function defined on a triangulated 2-dimensional or 3-dimensional domain \cite{EdelsbrunnerHarerZomorodian, BremerHamannEdelsbrunnerPasucci, GyulassyNatarajanPascucciHamann, GyulassyBremerHamannPascucci}. Furthermore, in \cite{GerberBremerPascucciWhitaker, GerberPotter} a Morse-Smale complex is statistically approximated from a high-dimensional scalar field. A difference between \cite{GerberBremerPascucciWhitaker, GerberPotter} and our setting is that their input is not a point cloud but instead a point cloud equipped with a function value at each point, thought of as a finite sampling from a scalar field. A second difference is that the scalar function in \cite{GerberBremerPascucciWhitaker, GerberPotter} is treated symmetrically, meaning that minima are as important as maxima, and the cells in the complex encode how one can descend from a maximum down to a minimum. Saddle points are not recovered, and any max-min pair is connected by at most one descending path. In our NEB approach, we equip the point cloud with a density function but we are not interested in the regions where this function is small or attains minima. Saddle points are of primary interest, as they correspond to high-density cells between maxima, and we are interested in detecting when there are multiple high-density paths between two maxima, for example as in Figure~\ref{fig: OpticalImagePatches_1cells} (center). For an application of the nudged elastic band method to the study of configuration spaces of hard disks, see \cite{HardDisks}.

\section{Background}\label{S: Background}

We briefly introduce three topics: CW complexes, Morse theory, and the nudged elastic band method.

\subsection{CW complexes}

\begin{figure}[h]
\centerline{\includegraphics[width=0.7\textwidth]{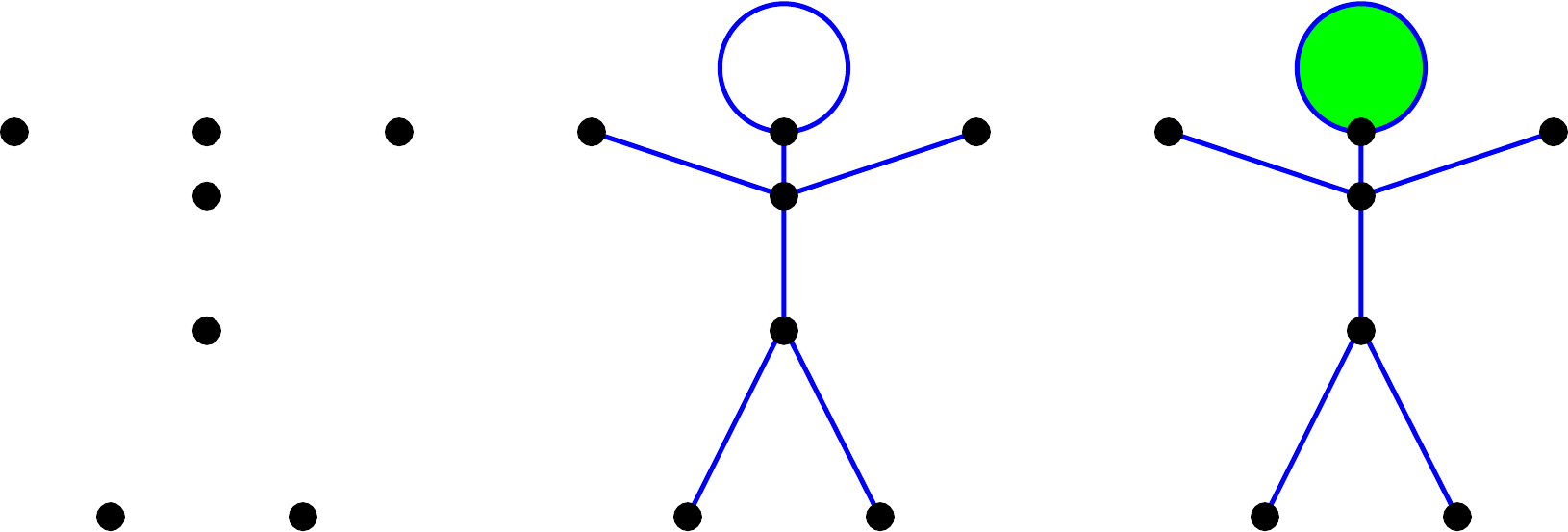}}
\caption{A stick figure represented as a CW complex containing eight 0-cells, eight 1-cells, and one 2-cell. The 0-skeleton is on the left, the 1-skeleton is in the center, and the full 2-skeleton is on the right.}
\label{fig: CWcomplex}
\end{figure}

A CW complex is a type of cell complex. For $k$ a nonnegative integer, a $k$-cell is the closed ball $\{y \in \R^k \;|\; \|y\| \leq 1\}$ of dimension $k$. So a 0-cell is a point, a 1-cell is a line segment, a 2-cell is a disk, etc. A CW complex $W$ is a topological space formed by the following inductive procedure. The 0-skeleton $W^{(0)}$ of $W$ is a set of 0-cells. The 1-skeleton $W^{(1)}$ is formed by gluing the endpoints of 1-cells to the 0-skeleton, and can be thought of as a graph. Inductively, we form the $k$-skeleton $W^{(k)}$ by gluing the boundaries of $k$-cells to the $(k-1)$-skeleton $W^{(k-1)}$. If $W$ is a finite-dimensional CW complex, then this process terminates and we have $W = W^{(k)}$ for some $k$. See Figure~\ref{fig: CWcomplex} for an example and \cite{Hatcher} for further details.

\subsection{Morse theory}

\begin{figure}[h]
\centerline{\includegraphics[width=0.9\textwidth]{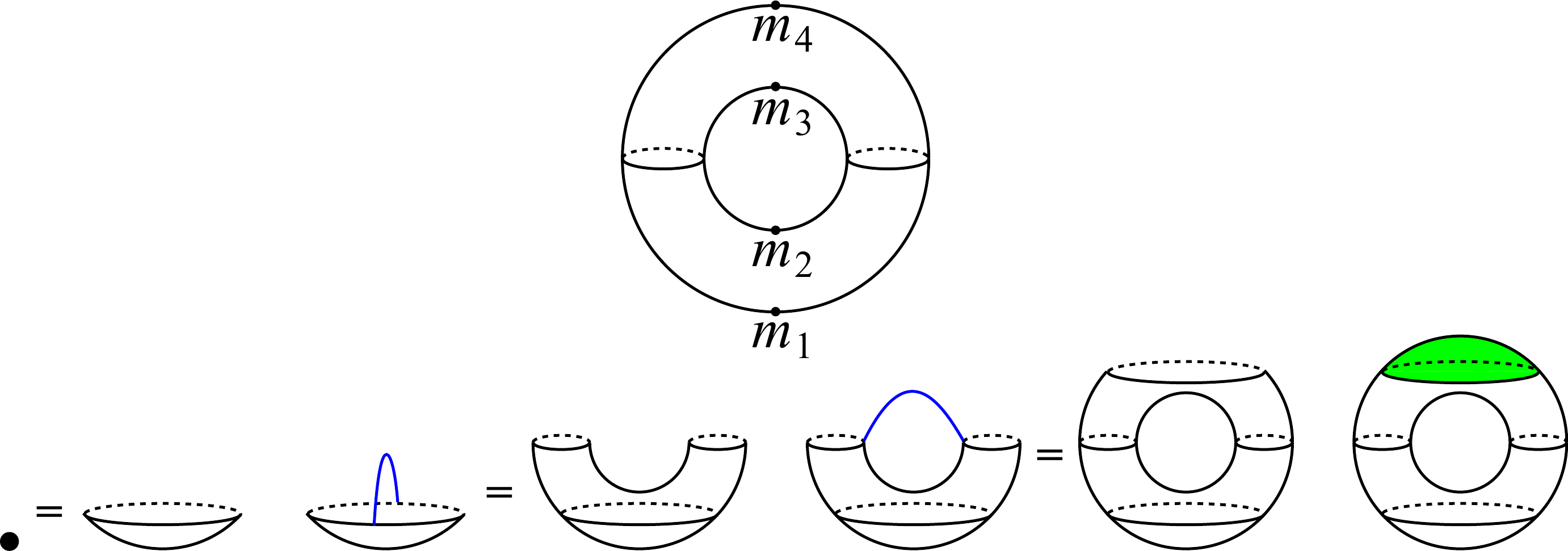}}
\caption{Morse theory example. ({\it Top}) Manifold $M$ is a torus and function $f \colon M \to \R$ is the height function. There are four critical points $m_1, \dots, m_4$ with indices 0, 1, 1, and 2 respectively. Let $a_0 < f(m_1) < a_1 \dots < f(m_4) < a_4$. ({\it Bottom}) Moving from left to right, we attach cells of dimension 0, 1, 1, and 2 in order to obtain $M_{a_i}$ from $M_{a_{i-1}}$.
}
\label{fig: MorseExample}
\end{figure}

The following introduction to Morse theory is informal; see \cite{Milnor} for a thorough treatment. Suppose $M$ is a compact manifold of dimension $d$ and Morse function $f \colon M \to \R$ is smooth with non-degenerate critical points $m_1, \dots, m_ k \in M$ satisfying
$$a_0 < f(m_1) < a_1 < f(m_2) < \dots < a_{k-1} < f(m_k) < a_k.$$
The index $\lambda_i$ of critical point $m_i$ is the number of linearly independent directions around $m_i$ in which $f$ decreases. So a minimum has index 0, a maximum has index $d$, and a saddle point has index between 1 and $d-1$. Let $M_a = f^{-1} \bigl( (-\infty, a] \bigr)$ be the sublevel set corresponding to $a \in \R$. Morse theory tells us that each $M_{a_i}$ is homotopy equivalent to a CW complex with one $\lambda_i$-cell for each critical point $m_i$. In particular, $M_{a_i}$ is homotopy equivalent to $M_{a_{i-1}}$ with a single $\lambda_i$-cell attached. For instance, $M_{a_1}$ is homotopy equivalent to a point and is obtained from $M_{a_0}$, the emptyset, by attaching a single 0-cell. Figure~\ref{fig: MorseExample} contains an example. 

Though Morse theory is traditionally stated in terms of sublevel sets, there is an equivalent formulation in terms of superlevel sets. The superlevel sets of $f$ correspond to the sublevel sets of $-f$, so $m_i$ is a critical point of $f$ with index $\lambda_i$ if and only if $m_i$ is a critical point of $-f$ with index $d-\lambda_i$. It follows that superlevel set $M^{a_{i-1}} = f^{-1} \bigl( [a_{i-1}, \infty) \bigr)$ is obtained from $M^{a_i}$ by attaching a single cell of dimension $d-\lambda_i$.

\subsection{Nudged elastic band}\label{SS: Nudged elastic band}

To find saddle points in a high-dimensional space we use the nudged elastic band method (NEB) from computational chemistry \cite{JonssonMillsJacobsen, HenkelmanJonsson, HenkelmanUberuagaJonsson, SheppardTerrellHenkelman}. Consider a chemical system, e.g.\ several molecules, whose space of states is parametrized by $\R^n$ and is equipped with a differentiable map $E \colon \R^n \to \R$ encoding the potential energy of the system at each state. The local minima of $E$ correspond to stable states, and chemists are interested in finding reaction paths between two stable states. A reaction path is a minimum energy path, whose points minimize the energy $E$ in all directions perpendicular to the path, and which passes through at least one saddle point of index one.

\begin{figure}[h]
\centerline{
\includegraphics[width=0.3\textwidth]{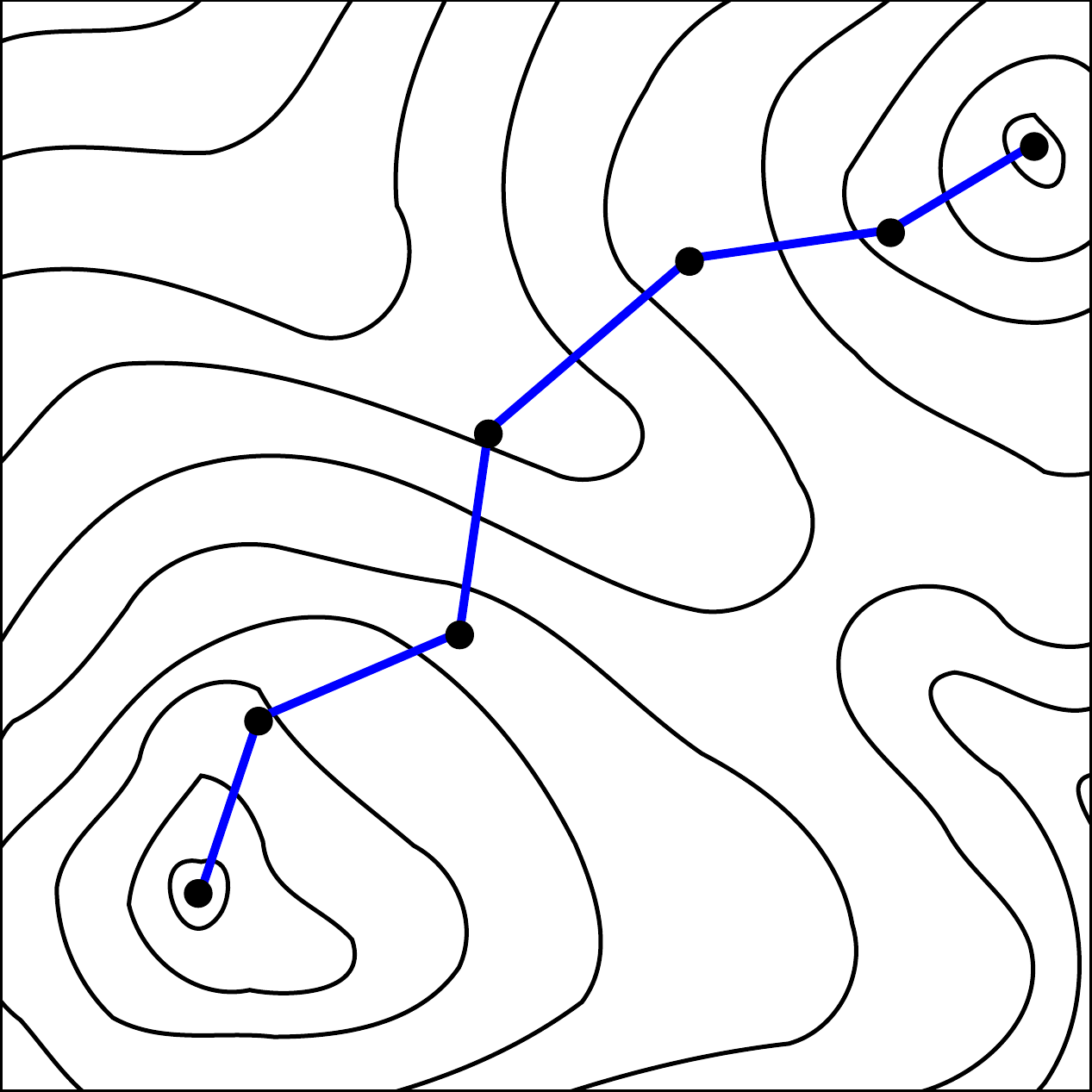}\hspace{15mm}\includegraphics[width=0.3\textwidth]{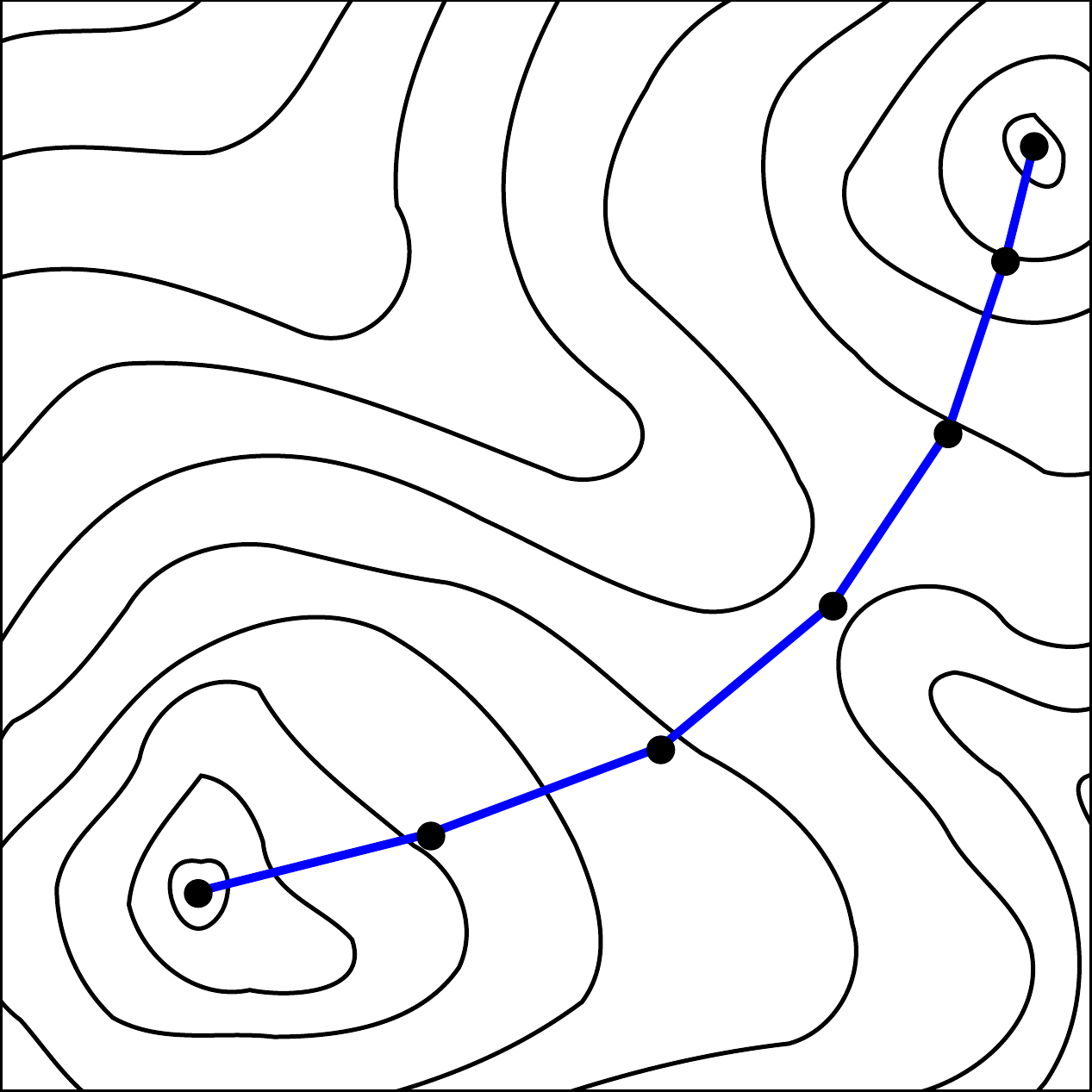}}
\caption{({\it Left}) An initial band connecting two local minima. ({\it Right}) A minimum energy path found using the nudged elastic band method.}
\label{fig: NudegedElasticBand_example}
\end{figure}
\FloatBarrier

The nudged elastic band method is used to find minimum energy paths. An initial piecewise linear path (called a band) connecting two local minima of energy $E$ is chosen, and forces move the band towards a minimum energy path. A gradient force moves each node in the band along the component of $-\nabla E$ perpendicular to the band, and a spring force ensures that adjacent edges of the band are of approximately equal length, so that the band does not tear apart. More explicitly, the gradient and spring forces appear as the first two terms on the right hand side of Equation (\ref{eqn: forces}) in Section \ref{SS: 1-cells}. Extra forces can be added to smoothen or dampen the motion of the band. Once the forces have stabilized, the convergent band is the approximation to the minimum energy path found by NEB.

\section{Method}\label{S: Method}

We suppose data set $X \subset \R^n$ is a finite sampling from an unknown probability density function $g \colon \R^n \to [0, \infty)$. We are interested in the regions of high density in the data set, and therefore would like to understand the profile of superlevel sets
$$Y^a = g^{-1} \bigl( [a, \infty) \bigr) = \{y \in \R^n \;|\; g(y) \geq a\}$$
containing all points in $\R^n$ with density at least $a \in [0, \infty)$. The superlevel sets $Y^a$ encode how the dense regions of data set $X$ are organized. 

We follow the ideas of Morse theory to build CW complex models $Z^a$ approximating the superlevel sets $Y^a$. First we use $X$ to build a differentiable density estimate $f \colon \R^n \to [0, \infty)$ approximating the unknown probability density function $g$. The 0-cells of our models will be local maxima of $f$. The 1-cells will be paths of high density between 0-cells, found using NEB. To find dense 2-cells with boundaries in the 1-skeleton we use an adaptation of NEB. Given a cell $e \subset \R^n$, let its density be $\inf \{f(y) \;|\; y \in e\}$. We define $Z^a$ to be the union of the cells with density greater than or equal to $a$.

Note that for $a \geq b$ we have the inclusion of superlevel sets $Z^a \subset Z^{b}$. A feature of our approach is that we compute the CW complex model $Z^a$ for many values of $a$ at once. This allows the user to observe how $Z^a$ changes as $a$ decreases. In accordance with the idea of topological persistence \cite{TopologicalPersistence, ComputingPersistent, ComputationalTopology}, knowing $Z^a$ over a range of values for $a$ gives a better picture of the data set than knowing $Z^a$ for any single choice of $a$. For example, in Section \ref{SS: Social network} we produce three different models---four vertices, a square loop, and a square disk---for a social network data set (Figure~\ref{fig: SocialNetwork_1cells}). Each model represents the data at a different density threshold, and together they provide a more complete understanding than any single CW complex model alone. In addition, our nested output can be used to measure the significance of topological features. Suppose a 1-cell with density $a$ appears in $Z^a$ to form a new loop, and suppose a 2-cell with density $b \leq a$ appears in $Z^b$ to fill this loop to a disk. Then $a-b$ measures how long this loop persists, and if $a-b$ is large then this loop is likely a significant feature of the data set, but if $a-b$ is small then the loop may be the product of noise. 

Several steps in our approach can be treated in a variety of ways. To name a few, we estimate density, find local maxima, generate random initial bands, simulate bands according to a formulation of NEB, and cluster convergent bands. We use simple methods to handle each of these steps, but we leave open the possibility of substituting more sophisticated methods. Our approach depends on the choice of several parameters, and the main parameter which we find necessary to tune is the standard deviation used to build a density estimator. We use a consistent choice for all other parameters across all of our data sets, which suggests the other parameters may not be as important to tune.

\subsection{Density estimation}

From data set $X \subset \R^n$, we build a differentiable density estimator $f \colon \R^n \to [0, \infty)$ approximating the unknown underlying density as follows. Let $\psi_{x, \sigma} \colon \R^n \to [0, \infty)$ be the probability density of a normal distribution centered at $x \in \R^n$ with standard deviation $\sigma > 0$. More explicitly, the $n \times n$ covariance matrix contains $\sigma^2$ along the diagonal and 0 elsewhere. We use the kernel estimator $f(y) = |X|^{-1} \sum_{x \in X}  \psi_{x, \sigma}(y)$. See \cite{Silverman} for other possible estimators, including different kernels or adaptive choices of standard deviation. Bandwidth estimators often aim to minimize the $L^2$ error, but this is not necessarily optimal for recovering the topology of sublevel sets. We regard standard deviation $\sigma$ as the main parameter the user must tune: one may increase $\sigma$ to smoothen $f$ or decrease $\sigma$ to expose local detail.

\subsection{0-cells}

To find the 0-cells of our model, we pick a uniformly random sample of points from $X$ and flow each along the gradient of $f$ towards a local maxima using the mean shift iterative procedure \cite{Fukunaga, Cheng}. If the size of $X$ is greater than 2,000 than we select 2,000 uniformly random points from $X$, and if the size of $X$ is at most 2,000 then we use each point of $X$ as a seed for mean shift. Let $y_0 \in X$ be one of the random initial points. We define a sequence $y_0, y_1, \dots$ by setting $y_{i+1} = m(y_i)$, where the mean shift function $m \colon  \R^n \to \R^n$ is defined by 
$$m(y) = \frac{\sum_{x \in X} \psi_{x, \sigma}(y) x}{\sum_{x \in X} \psi_{x, \sigma}(y)}\ .$$
The vector $m(y) - y$ is proportional to the normalized gradient $\nabla f(y) / f(y)$, and the $y_i$ converge to a local maxima of $f$. It is necessary to cluster the convergent points to identify which represent the same 0-cell, and we use single-linkage clustering. Since we are interested in CW complex models with more than one 0-cell, we select the clustering threshold parameter to lie in longest interval where a constant number of clusters greater than one is obtained. For the social network, optical image, range image, and optical flow examples in Sections~\ref{SS: Social network}--\ref{SS: Optical flow patches}, the longest interval with a constant number of clusters constitutes 66.0\%, 97.8\%, 99.6\%, and 90.6\% of the total parameter range with more than one cluster, respectively. We select the densest member from each cluster as a 0-cell in our model.

\subsection{1-cells}\label{SS: 1-cells}

To find the 1-cells of our model we use NEB, which we now describe in our data analysis setting. Our formulation is similar to \cite{JonssonMillsJacobsen, HenkelmanJonsson, HenkelmanUberuagaJonsson, SheppardTerrellHenkelman}. A piecewise linear band is given by a list of nodes $v_1, v_2, \dots, v_N$, with endpoints $v_1$ and $v_N$ in our set of 0-cells. Forces act on the intermediate nodes $v_i$ with $1 < i < N$ while the endpoints remain fixed. The first task is to approximate a unit tangent vector $\tau_i$ at each intermediate node. Define $u_i^+ = v_{i+1} - v_i$ and $u_i^- = v_i - v_{i-1}$. We use a na\"ive tangent estimate $\tau_i = (u_i^+ + u_i^-) / \|u_i^+ + u_i^-\|$ given by averaging, though more elaborate tangent estimators may be used \cite{HenkelmanJonsson}.

At each intermediate node $v_i$ we define a total force
\begin{equation}\label{eqn: forces}
F_i\ =\ c \nabla f(v_i)|_\perp\ +\ (\|u_i^+\| - \|u_i^-\|) \tau_i\ +\ smoothing.
\end{equation}
The expression $\nabla f(v_i)|_\perp$ is the component of $\nabla f(v_i)$ perpendicular to the tangent $\tau_i$, and is called the gradient force. The gradient constant $c$ adjusts the strength of the gradient force. To normalize with respect to the maximum gradient of the normal distribution, we set
$$c = \bigl( \sup_{y \in \R^n} \| \nabla \psi_{\vec{0}, \sigma}(y) \| \bigr)^{-1} = \bigl( \sigma \sqrt{2 \pi} \bigr)^n \sqrt{e}.$$
The term $(\|u_i^+\| - \|u_i^-\|) \tau_i$ in (\ref{eqn: forces}) is the spring force. This name is slightly misleading, as the spring force neither enforces a natural spring length on each edge nor minimizes the length of each edge. Instead, the spring force aims to equate the lengths of adjacent edges in the band.

The {\it smoothing} term in (\ref{eqn: forces}) is a force added to prevent kinks, which hinder convergence, from forming in the band. Let $\theta_i$ be the angle between vectors $u_i^+$ and $u_i^-$. Typically these vectors are close to parallel and $\theta_i$ is close to zero. Let $0 \leq \alpha < \beta \leq \pi$ be fixed angles, and let $h_{\alpha, \beta} \colon [0, \pi] \to [0, 1]$ be a function which is zero for $x \leq \alpha$, one for $x \geq \beta$, and increasing continuously from zero to one as $x$ increases from $\alpha $ to $\beta$. We define our smoothing force to be $h_{\alpha, \beta}(\theta_i)(u_i^+ - u_i^-)$, which moves $v_i$ in order to decrease $\theta_i$ whenever $\theta_i > \alpha$. In this work we set
$$
h_{\alpha, \beta}(x)=
\begin{cases}
0 & \textrm{if } x \leq \alpha \\
\bigl( 1 - \cos(\pi \frac{x - \alpha}{\beta - \alpha}) \bigr) / 2 & \textrm{if } \alpha < x < \beta \\
1& \textrm{if } x \geq \beta
\end{cases}
$$
with $\alpha = \pi / 6$ and $ \beta = \pi / 2$. We obtain the same CW complex models for the social network, optical image, range image, and optical flow examples in Sections~\ref{SS: Social network}--\ref{SS: Optical flow patches} using the smoothing force in \cite{JonssonMillsJacobsen}, although occasionally kinks appear in the bands prior to convergence.

\begin{figure}[h]
\centerline{\includegraphics[width=.45\textwidth]{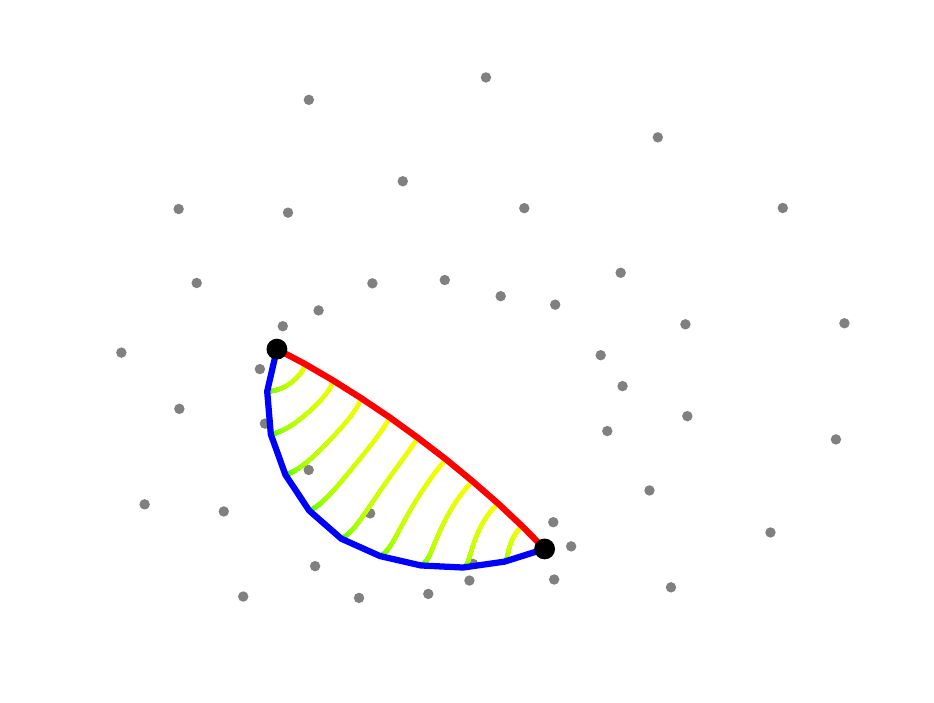}}
\caption{A sample 1-cell trial on a small data set. The initial band is in red and the convergent 1-cell is in blue. The lines which fade from yellow to green trace the paths of the intermediate nodes.}
\label{fig: sampleTrial}
\end{figure}

Evolving the band amounts to numerically solving the system of first order differential equations $v_i' = F_i$; see Figure~\ref{fig: sampleTrial}. In the chemistry setting it is appropriate to set acceleration proportional to the gradient of the potential energy. In our setting we do not view $\nabla f$ as a force in the literal sense but only as an indication of which direction to move in order to maximize $f$, and hence we use the first derivative $v_i'$ rather than the second derivative. Nevertheless, we have also had success with the second order equation and do not dismiss its use.

Between two 0-cells $p$ and $q$  in our model there may be no 1-cells, a single 1-cell, or multiple 1-cells. To find the 1-cells we generate a sample of initial bands joining $p$ and $q$, as described in Appendix \ref{A: Initial bands}. We evolve each band until it converges and discard non-convergent bands. We also discard convergent bands which pass too close to any other 0-cell $r \neq p, q$, as such a band should instead appear as the concatenation of two others. To identify which convergent bands represent the same 1-cell, we perform single-linkage clustering on the bands with threshold parameter equal to 10\% of the diameter of the data set. The social network, optical image, range image, and optical flow examples in Sections~\ref{SS: Social network}--\ref{SS: Optical flow patches} are quite stable with respect to this choice --- we obtain identical results for any percentage of the diameter between 1\% and 45\%. Our metric $d_{p, q, N}$ on bands of $N$ nodes starting at $p$ and ending at $q$ assigns to the two bands $v_1, \dots, v_N$ and $\tilde{v}_1, \dots, \tilde{v}_N$ the distance $d_{p, q, N}(\{v_i\}, \{\tilde{v}_i\}) = (N-2)^{-1} \sum_{i=2}^{N-1} d(v_i, \tilde{v}_i)$. From each cluster we select the band with the highest density to represent the corresponding 1-cell. We estimate a band's density as $\min \{f(v_1), \dots, f(v_N)\}$, though one could obtain a more accurate estimate by subdividing the band further or by using the climbing image method of \cite{HenkelmanUberuagaJonsson}.

We make the following consistent parameter choices for all of the data sets. We include $N=11$ nodes in the bands and say that a band has converged when $(N-2)^{-1} \sum_{i=2}^{N-1} \| v_i' \| < 10^{-4}$. Increasing $N$ or decreasing the convergence threshold improves the maximum density path approximation at the expense of computation time, but does not alter our CW complex models.

\subsection{Higher-dimensional cells} 

In order to search for higher-dimensional cells, one can imagine adapting NEB. In Appendix \ref{A: Higher-dimensional cells} we describe a na\"ive approach along with its weaknesses. We are interested in better methods, but further work is needed.

\section{Results}\label{S: Results}

We apply our method to the four data sets in Table \ref{table: data sets}. From the social network, optical image, and range image data sets we extract small cell complexes which efficiently model the dense regions of the data. Moreover, we identify new topological structure in the data set of optical flow patches. See Appendix \ref{A: Additional data set information} for further information on the data sets. Our code is available along with several of the data sets at the following webpage: \url{http://code.google.com/p/neb-tda}.

\begin{table}[h]
\begin{threeparttable}
\caption{Data set information}\label{table: data sets}
\begin{tabular}{@{\vrule height 10.5pt depth4pt  width0pt}lcccc}
\toprule
& social network & optical image & range image & optical flow \\
\midrule
size of data set $X \subset \R^n$ & 1,127 & 15,000 & 15,000 & 15,000 \\
dimension $n$  & 5 & 8 & 24 & 16 \\
standard deviation\tnote{1}\ \ $\sigma$ & 0.45 & 0.20 & 0.35 & 0.30 \\
\bottomrule
\end{tabular}
\begin{tablenotes}[para]
\item[1] This is not a property of the data set, but instead the standard deviation we use to estimate density.
\end{tablenotes}
\end{threeparttable}
\end{table}

\subsection{Social network}\label{SS: Social network}

The National Longitudinal Study of Adolescent Health is a school-based study of American youth. In 1994-95 a sample of high schools and middle schools was chosen, and when possible each high school was paired with a sister middle school from the same community. The students from each school were asked to list up to five of their closest female friends who attend their school or their sister school, and likewise for their male friends.

In \cite{Moody}, Moody creates network graphs from the survey data. Each student is represented by a vertex, and the edge between two students exists if each student listed the other as a friend. We analyze the graph for the students from ``Countryside High School,'' a pseudonym used by Moody, and its sister middle school. Of the 1,147 students from this community who participated in the survey, 20 students shared no connections with any others and were removed from the graph.

\begin{figure}[h]
\centerline{\includegraphics[width=0.9\textwidth]{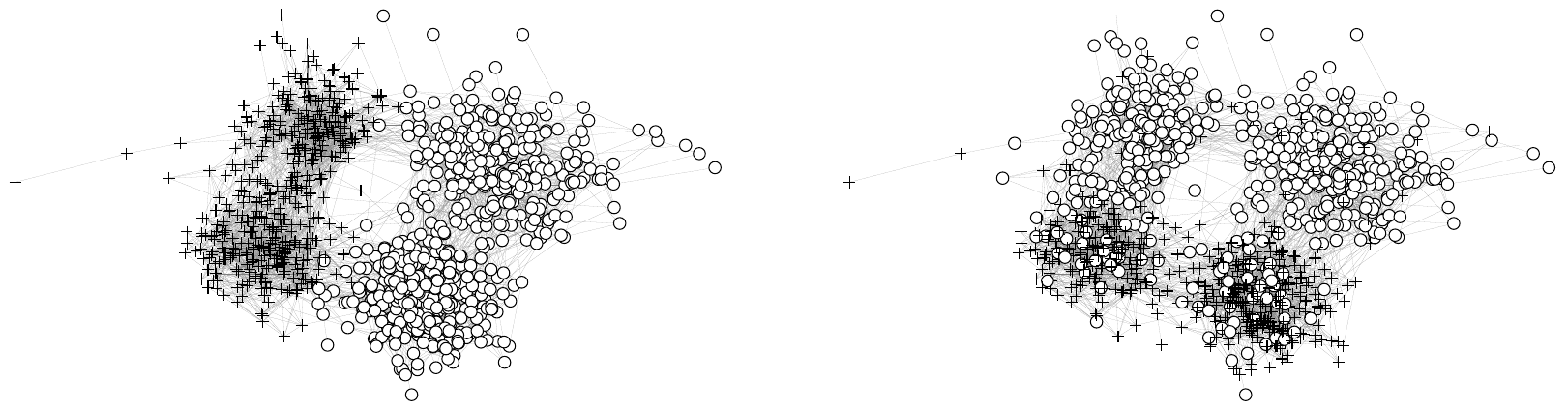}}
\caption{Social network for ``Countryside High School'' and its sister middle school. ({\it Left}) Cross vertices are students from the middle school and circle vertices are students from the high school. ({\it Right}) Cross vertices are white students and circle vertices are nonwhite students. A handful of students are without race data and their vertices are left unmarked.}
\label{fig: SocialNetwork_scode_race}
\end{figure}

Figure~\ref{fig: SocialNetwork_scode_race} illustrates the following interesting structure in this graph. Most of the students can be placed into one of four groups, containing a majority of white high schoolers, nonwhite high schoolers, white middle schoolers, or nonwhite middle schoolers. There are many friendships between students in the same group. In addition, there are a significant number of friendships between groups of the same school or race category. However, their are very few friendships between groups with neither category in common.

We illustrate how our method can recover relevant structure in this graph. The shortest path distance between vertices $v$ and $w$ is the fewest number of edges one must cross to travel from $v$ to $w$. To represent the graph as a point cloud we use stress majorization, an optimization strategy in multidimensional scaling \cite{CoxCox}, to embed the vertices of the graph as a data set $X \subset \R^5$ in a manner distorting the shortest path metric as little as possible. The choice of embedding dimension is not critical in this example --- we have tested our method on embeddings in $\R^n$ for $2 \leq n \leq 8$ and obtained the same CW complex models, where higher choices of $n$ allow for smaller distortion of the shortest path metric. Each point $x \in X$ corresponds to a student, and we regard the resulting set $X$ as our input data set.

\begin{figure}[h]
\centerline{\includegraphics[width=1.0\textwidth]{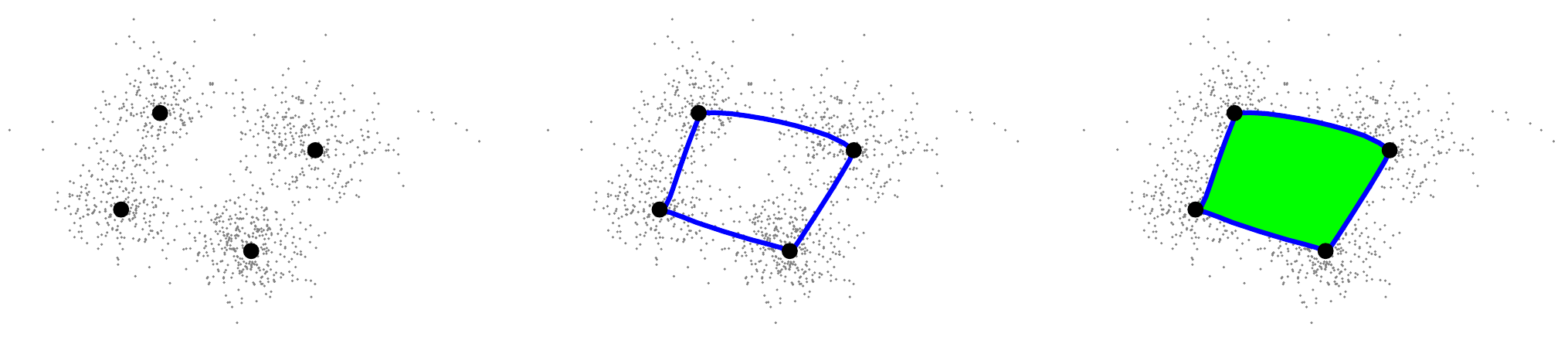}}
\caption{Social network data, projected to a plane using principal component analysis. Our CW complex model $Z^a$ grows as the density threshold $a$ decreases. From left to right, we have four 0-cells, a square loop, and a disk.}
\label{fig: SocialNetwork_1cells}
\end{figure}

Our method finds four 0-cells, which correspond to the four groups of students. We find four 1-cells which form a square loop. The 0-cells have densities in $[1.69 \cdot 10^{-2},\ 2.02 \cdot 10^{-2}]$ and the 1-cells in $[1.07 \cdot 10^{-2},\ 1.54 \cdot 10^{-2}]$.\footnote{Occasionally we also find a 0-cell with density below $6 \cdot 10^{-4}$, but due to its low density value it does not affect our analysis.} We find a 2-cell filling the loop with density $0.79 \cdot 10^{-2}$. Recall CW complex $Z^a$ is defined to be the union of the cells with density greater than or equal to $a$, and see Figure~\ref{fig: SocialNetwork_1cells} for three such CW complex models. For $a \in (1.54 \cdot 10^{-2},\ 1.69 \cdot 10^{-2})$ the model $Z^a$ consists of four 0-cells. Hence we recover the groupings of students based on school and race. For $a \in (0.79 \cdot 10^{-2},\ 1.07 \cdot 10^{-2})$ the model $Z^a$ is a square. The square reveals that groups sharing either the school or race category are more closely linked than groups sharing neither. This suggests that when making friends it is more difficult to cross two cultural barriers than one. For  $a < 0.79 \cdot 10^{-2}$ the model $Z^a$ fills to a disk, which is an appropriate representation of the data at a sufficiently coarse scale. Note that this increasing sequence of CW complex models provides a better understanding of the data set than any single model alone.

\FloatBarrier

\subsection{Optical image patches}\label{SS: Optical image patches}

\begin{figure}[h] 
\centerline{\includegraphics[width=1.0\textwidth]{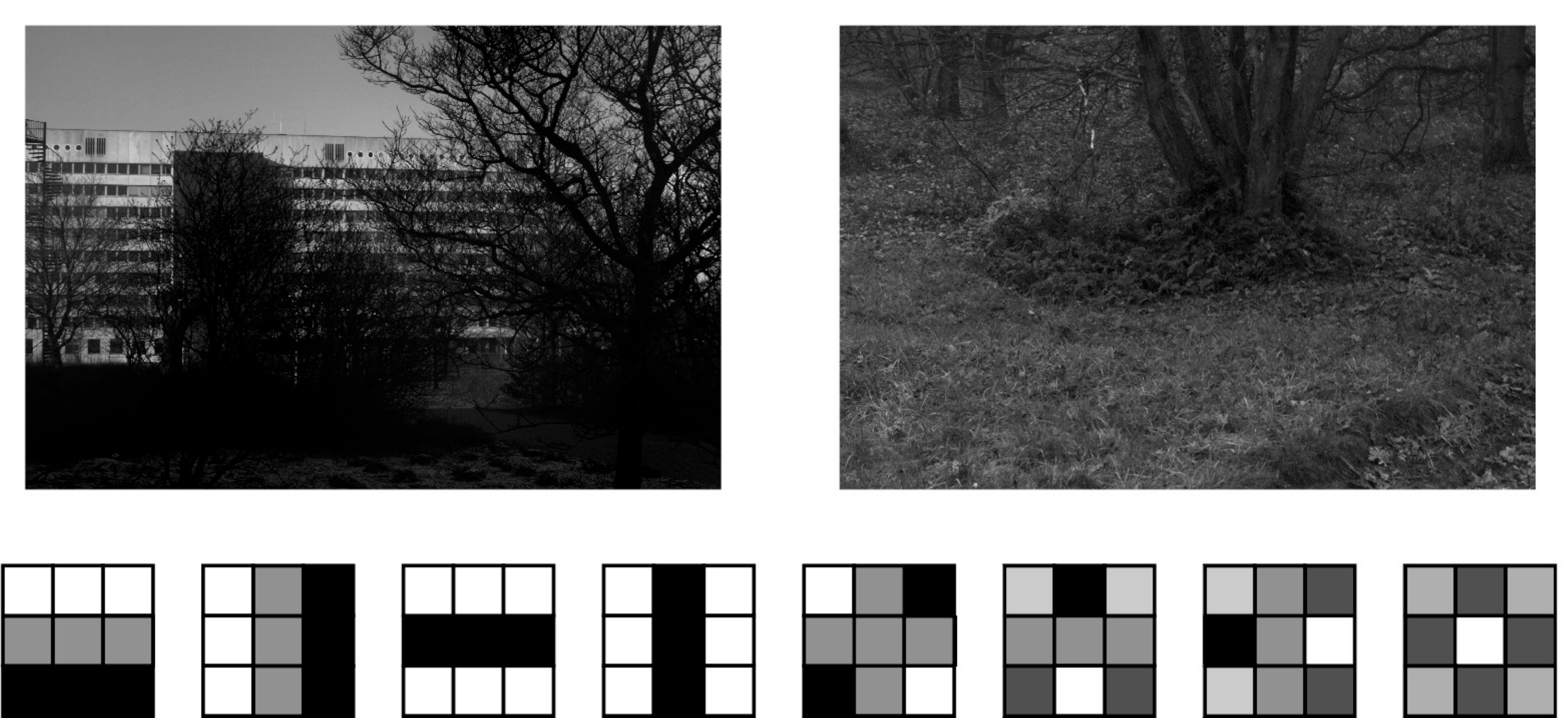}}
\caption{({\it Top}) Sample optical images from the van Hateren database. ({\it Bottom}) Vectors $e_1, e_2, \dots, e_8$ form the $3 \times 3$ DCT basis.}
\label{fig: opticalImages}
\end{figure}

The optical image database collected by van Hateren and van der Schaaf in \cite{VanHaterenVanDerSchaaf} contains a variety of indoor and outdoor scenes. From this database, Lee et al.\ in \cite{LeePedersenMumford} select a large random sample of $3 \times 3$ patches, each thought of as a point in $\R^9$. The coordinates are the logarithms of grayscale pixel values, and they use the logarithm so that the relative reflectance of an object is invariant under changes in lighting intensity. Lee et al.\ define a norm measuring the contrast of a patch and select the high-contrast patches. They normalize each patch by subtracting from each pixel the average of the nine coordinates, and by dividing by the contrast norm. They change to the Discrete Cosine Transform (DCT) basis $\{e_1, e_2, \dots, e_8\}$, which maps the patches to the unit sphere $S^7 \subset \R^8$. Note in Figure~\ref{fig: opticalImages} that $e_1$ and $e_2$ are horizontal and vertical linear gradients while $e_3$ and $e_4$ are horizontal and vertical quadratic gradients. Let $\cM$ be the resulting set of high-contrast, normalized, $3 \times 3$ patches. Each point of $S^7$ is close to some point of $\cM$, but the density of $\cM$ varies widely.

\begin{figure}[h]
\centerline{\includegraphics[width=0.45\textwidth]{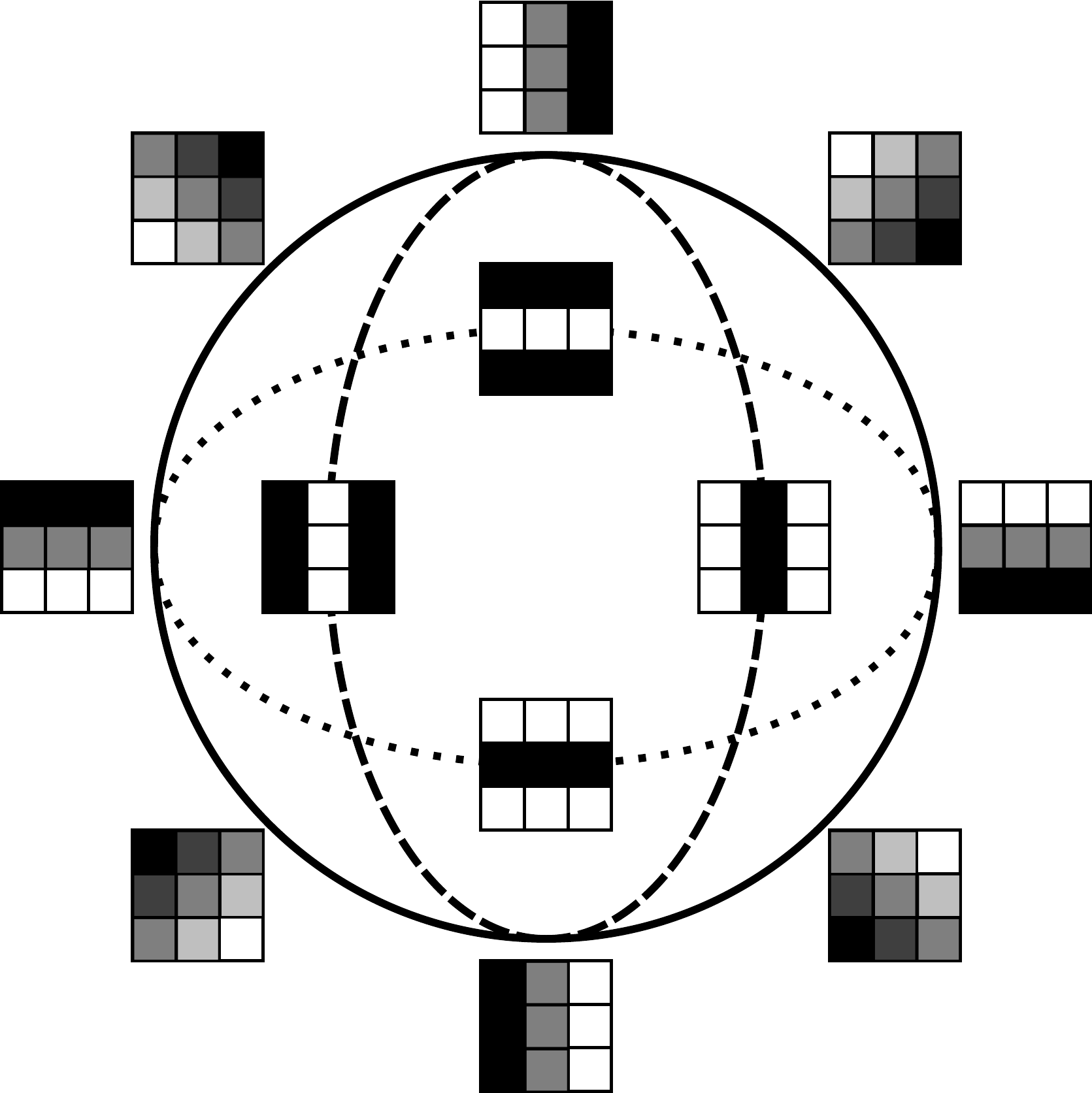}}
\caption{Three circle model. The solid outer circle is the primary circle and contains linear gradients. The dotted and dashed inner circles are the horizontal and vertical secondary circles which contain quadratic gradients. Each secondary circle intersects the primary circle twice, but the secondary circles do not intersect each other.}
\label{fig: threeCircle}
\end{figure}

Carlsson et al.\ in \cite{KleinBottle} use persistent homology to study $\cM$. Using a family of density estimators they select a family of dense core subsets from $\cM$. They apply persistent homology, and for a global estimate of density their core subset has $\rank(H_0) = \rank(H_1) = 1$, where $H_i$ denotes the $i$-th homology group. This core subset lies near the primary circle $\{\alpha e_1 + \beta e_2 \;|\; (\alpha, \beta) \in S^1\}$ containing linear gradient patches at all angles. With a more local estimate of density, the core subset has $\rank(H_0) = 1$ and $\rank(H_1) = 5$. Carlsson et al.\ identify a three circle model matching this homology profile and the data; see Figure~\ref{fig: threeCircle}. This is a supervised step. In addition to the primary circle, the three circle model contains two secondary circles, $\{\alpha e_1 + \beta e_3 \;|\; (\alpha, \beta)\in S^1\}$ and $\{\alpha e_2 + \beta e_4 \;|\; (\alpha, \beta) \in S^1\}$, which include quadratic gradients in the horizontal or vertical direction. The primary circle reflects nature's preference for linear gradients in all directions, and the secondary circles reflect nature's preference for the horizontal and vertical directions.

\begin{figure}[h]
\centerline{\includegraphics[width=0.75\textwidth]{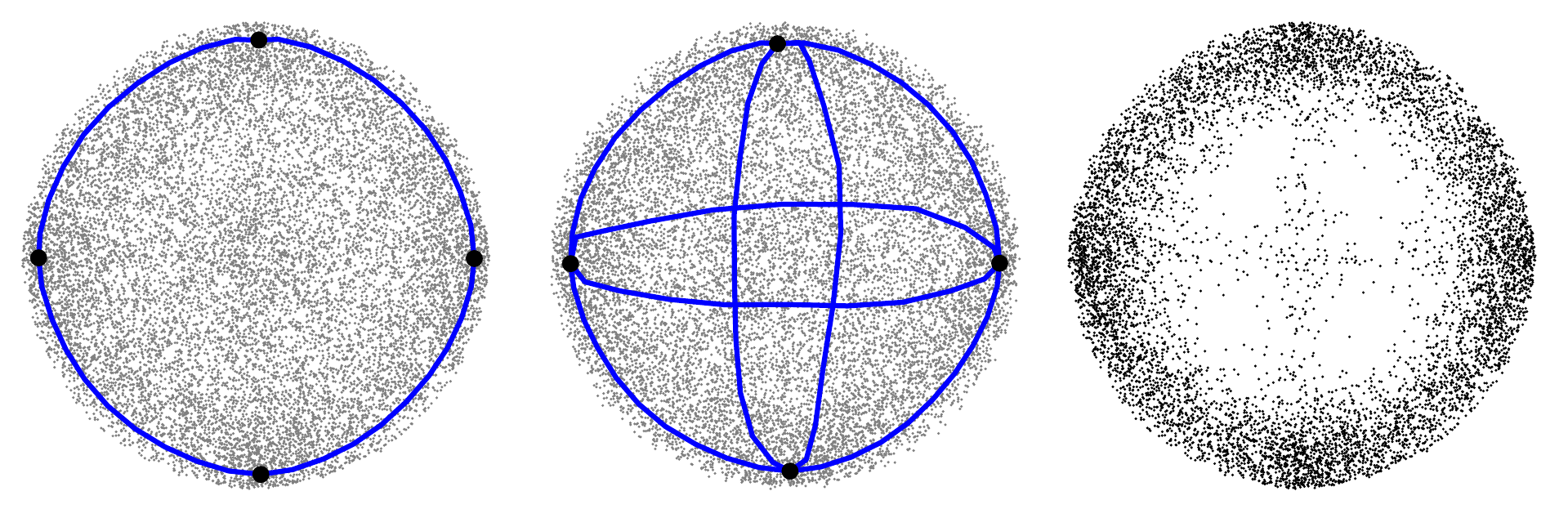}}
\caption{({\it Left}) Optical image patches $X \subset \R^8$ and primary circle complex $Z^a$ for $a \in (0.38, 1.17)$, projected to the $e_1 e_2$ plane. ({\it Center}) Three circle model $Z^a$ for $a < 0.33$, projected to the plane spanned by $e_1 + \frac{1}{4} e_4$ and $e_2 + \frac{1}{4} e_3$ so that the secondary circles are visible. ({\it Right}) The densest 8,500 points of $X$, projected to the $e_1 e_2$ plane. The primary circle appears clearly and the projections of the secondary circle patches form a faint cross.}
\label{fig: OpticalImagePatches_1cells}
\end{figure}

For our analysis we select a random input data set $X \subset \cM$ of size 15,000. We find four 0-cells located near the four most common patches $\pm e_1$ and $\pm e_2$. Between each of the four 0-cell pairs $\{e_1, e_2\}$, $\{e_2, -e_1\}$, $\{-e_1, -e_2\}$, and $\{-e_2, e_1\}$ we find a quarter-circular 1-cell. Together these form the primary circle. Between the pair $\{e_1, -e_1\}$ we find two semicircular 1-cells forming the horizontal secondary circle, and between $\{e_2, -e_2\}$ we find two semicircular 1-cells forming the vertical secondary circle. The 0-cells have densities in  $[2.11, 2.36]$,\footnote{Occasionally we also find a 0-cell with density below $0.08$, but due to its low density value it does not affect our analysis.} the primary circle 1-cells in $[1.17, 1.28]$, and the secondary circle 1-cells in $[0.33, 0.38]$. For $a \in (1.28, 2.11)$ our model $Z^a$ is the four most common patches, for $a \in (0.38, 1.17)$ it is the primary circle, and for $a < 0.33$ it is the three circle model (Figure~\ref{fig: OpticalImagePatches_1cells}). Indeed, in an unsupervised fashion we obtain the minimal cellular decomposition of the three circle model that is geometrically accurate.

Carlsson et al.\ discover a 2-dimensional Klein bottle surface that contains the three circle model as its backbone and which is a good example of how data set models can not only improve qualitative understanding but also lead to applications. As a low-dimensional manifold the Klein bottle can be used in image compression schemes \cite{KleinBottle}, and in addition the Klein bottle model has been used to identify and analyze optical image textures \cite{PereaCarlsson}. In order to identify the Klein bottle using persistent homology, Carlsson et al.\ add a set $Q \subset \cM$ of optical image patches to their data set, where ``the size of Q can be thought of as a measure of the failure of the density function to cut out the 2-dimensional manifold" \cite{KleinBottle}. As this is a first exploration of nudged elastic band in topological data analysis, we do not add any such points $Q$ in order to try to corroborate the Klein bottle model.
\FloatBarrier

\subsection{Range image patches}

A range image is captured by a laser scanner, and each pixel in a range image stores the distance between the laser scanner and the nearest object in the corresponding direction. In \cite{LeePedersenMumford}, Lee et al.\ select a large random sample of log-valued, high-contrast, normalized, $3 \times 3$ range image patches from the Brown database, which contains a variety of indoor and outdoor scenes. They observe that the patches cluster near binary patches, where the binary values correspond to foreground and background.

\begin{figure}[h]
\centerline{\includegraphics[width=.7\textwidth]{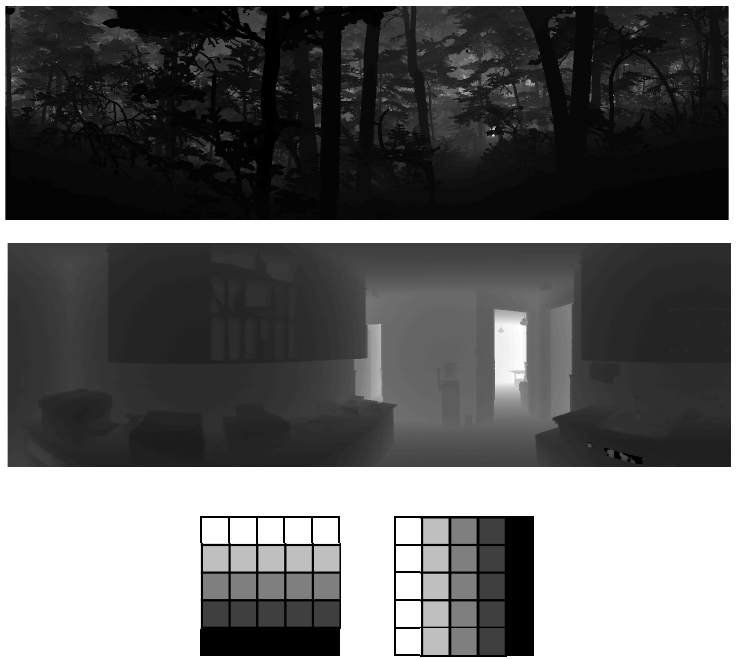}}
\caption{({\it Top}) Sample range images from the Brown database. ({\it Bottom}) Vectors $e_1$ and $e_5$ from the $5 \times 5$ DCT basis are horizontal and vertical linear gradients.}
\label{fig: rangeImages}
\end{figure}

Though the largest clusters are arranged in the shape of a circle, the $3 \times 3$ binary patches are too coarse to fill the full circle. In \cite{Range}, Adams and Carlsson consider $5 \times 5$ patches, preprocessed in manner similar to \cite{LeePedersenMumford}. They obtain a large sample $\cR$ of high-contrast, normalized, $5 \times 5$ range image patches. The $5 \times 5$ DCT basis now contains 24 vectors $\{e_1, e_2, \dots, e_{24}\}$, where $e_1$ and $e_5$ are horizontal and vertical linear gradients (Figure~\ref{fig: rangeImages}). With persistent homology they find that a dense core subset of $\cR$ has $\rank(H_1) = 1$, and they propose the range primary circle model $\{\alpha e_1 + \beta e_5 \;|\; (\alpha, \beta) \in S^1\}$ to match this homology.

\begin{figure}[h]
\centerline{\includegraphics[width=0.7\textwidth]{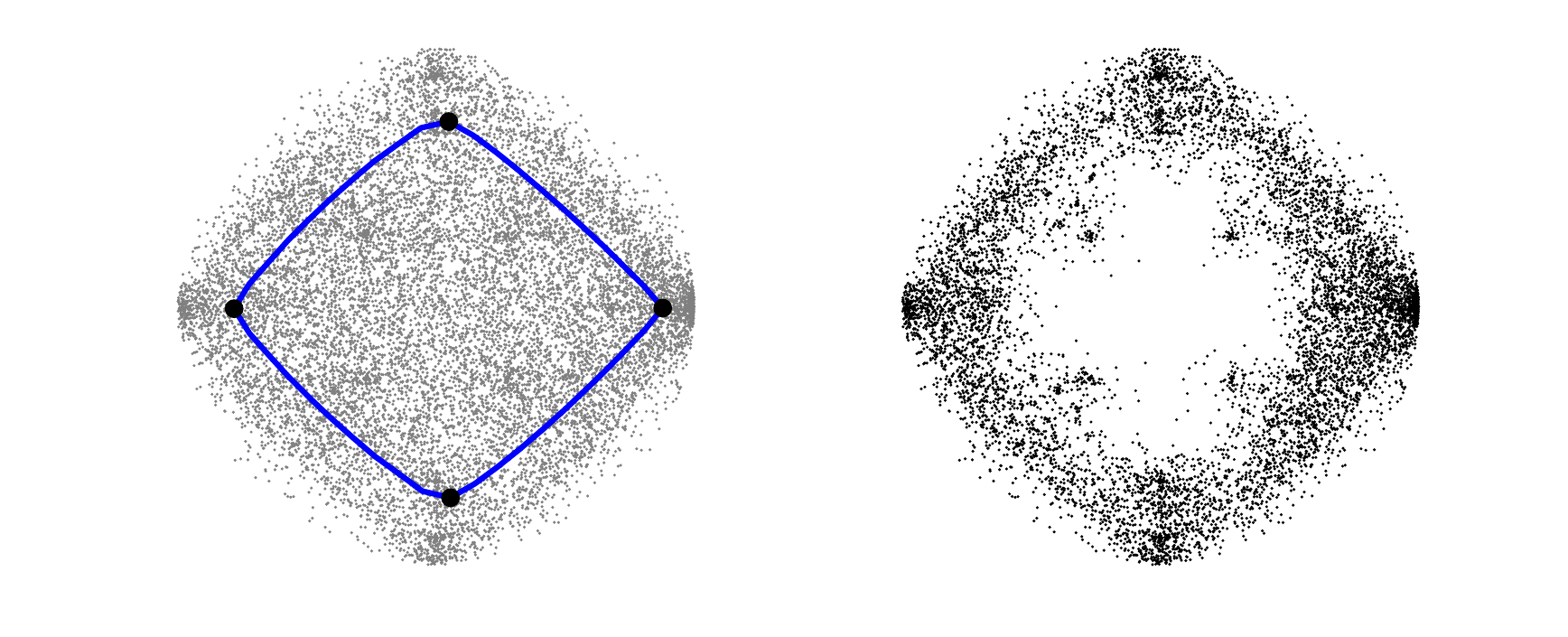}}
\caption{({\it Left}) Range image patches $X \subset \R^{24}$ and range primary circle $Z^a$ for $a \in (0.39, 0.51)$, projected to the $e_1 e_5$ plane. ({\it Right}) The densest 9,000 points of $X$.}
\label{fig: RangeImagePatches_1cells}
\end{figure}

For our Morse theory approach, we pick a random subset $X \subset \cR$ of size 15,000. We find four 0-cells near $\pm e_1$ and $\pm e_5$. Three of these 0-cells have densities in $[0.59, 0.75]$ while the 0-cell near $e_1$ has density 1.47. This reflects the fact that many range patches are shots of the ground and hence are near the horizontal linear gradient given by $e_1$. We find four quarter-circular 1-cells forming a loop with densities in $[0.51, 0.64]$, and a 2-cell filling the loop with density $0.39$.
For $a \in (0.75, 1.47)$ our model $Z^a$ is the single 0-cell near DCT basis vector $e_1$, and for $a \in (0.39, 0.51)$ our model $Z^a$ is the range primary circle (Figure~\ref{fig: RangeImagePatches_1cells}). 

\FloatBarrier

\subsection{Optical flow patches}\label{SS: Optical flow patches}

In this example we identify new topological structure in an optical flow data set. A video records a sequence of images, and its optical flow is the apparent motion in the sequence of images. At each frame, the optical flow is represented by a vector field with one vector per pixel that points to where that pixel appears to move for the subsequent frame. No instrument measures optical flow directly, and estimating optical flow from a video is an ill-posed problem. However, Roth and Black in \cite{RothBlack} create a database of ground-truth optical flow by pairing range images with camera motions and by calculating the produced optical flow. The range images are from the Brown database, and the camera motions are retrieved from a database of videos from hand-held or car-mounted cameras.

\begin{figure}[h]
\centerline{\includegraphics[width=.5\textwidth]{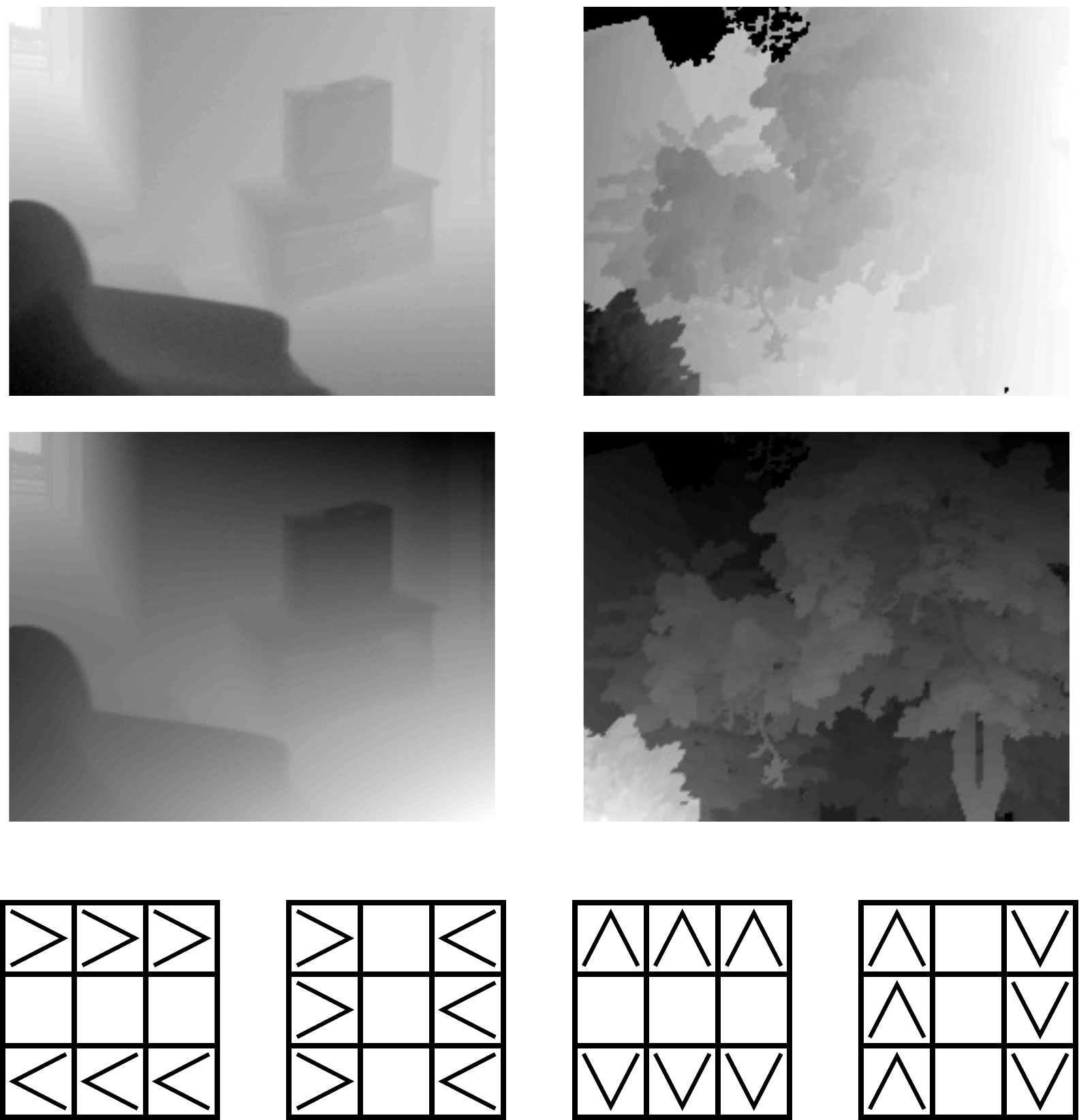}}
\caption{({\it Top}) Two sample optical flows from the Roth and Black database. Horizontal components are on the top and vertical components are on the bottom. White corresponds to flow in the positive direction ($+x$ or $+y$) and black corresponds to the negative direction. ({\it Bottom}) Vector fields $e_1^u$, $e_2^u$, $e_1^v$, and $e_2^v$ from the $3 \times 3$ optical flow DCT basis. Compare with the patches $e_1$ and $e_2$ in Figure~\protect\ref{fig: opticalImages}.}
\label{fig: opticalFlowImages}
\end{figure}

With preprocessing steps analogous to \cite{LeePedersenMumford} we create a large sample $\cF$ of high-contrast, normalized, $3 \times 3$ optical flow patches. We change to the DCT basis $\{e_1^u, \dots, e_8^u, e_1^v, \dots, e_8^v\}$, where the superscript $u$ denotes flow in the horizontal direction and $v$ denotes the vertical direction; see Figure~\ref{fig: opticalFlowImages}. We select $X$ to be a random subset of $\cF$ of size 15,000.

\begin{figure}[h]
\centerline{\includegraphics[width=0.7\textwidth]{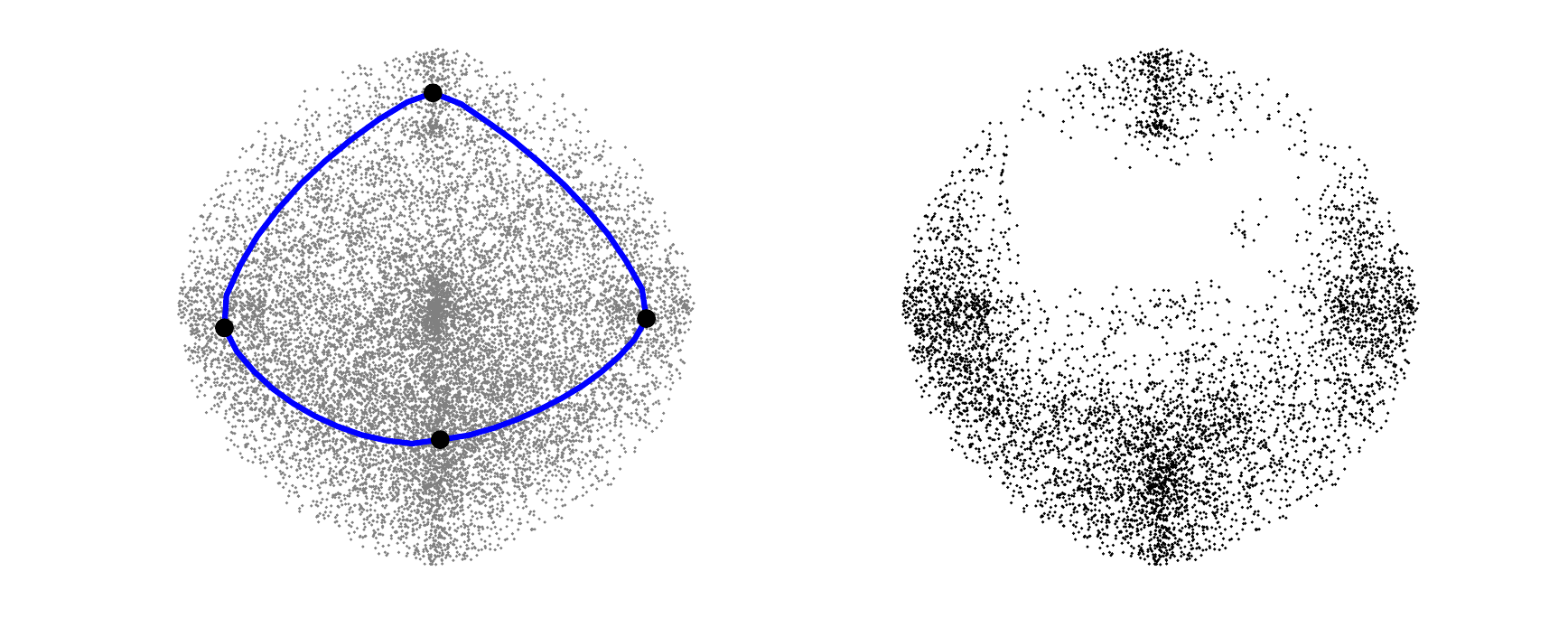}}
\caption{({\it Left}) Optical flow patches $X \subset \R^{16}$ and horizontal flow circle $Z^a$ for $a \in (0.39, 0.71)$, projected to the $e_1^u e_2^u$ plane. ({\it Right}) The densest 6,000 points of $X$.}
\label{fig: OpticalFlowPatches_1cells}
\end{figure}

We apply our method and find four 0-cells near $\pm e_1^u$ and $\pm e_2^u$ and four quarter-circular 1-cells, all with densities in $[0.71, 2.79]$. These cells form a loop. We find a 2-cell filling the loop with density 0.39, and hence for $a \in (0.39, 0.71)$ our model $Z^a$ recovers the horizontal flow circle near $\{\alpha e_1^u + \beta e_2^u \;|\; (\alpha, \beta) \in S^1\}$ (Figure~\ref{fig: OpticalFlowPatches_1cells}). When applying a camera motion to a range image, the apparent motion of the foreground is faster than that of the background. In particular, applying rightward horizontal camera motion to the linear $3 \times 3$ range patch $\alpha e_1 + \beta e_2$ (with $e_1$ and $e_2$ from the $3 \times 3$ DCT basis) produces the optical flow patch $\alpha e_1^u + \beta e_2^u$, and applying leftward horizontal camera motion to this same range patch produces the optical flow patch $-\alpha e_1^u - \beta e_2^u$. Hence the horizontal flow circle can be obtained by applying horizontal camera motion to range patches from the $3 \times 3$ analogue of the range primary circle. Also, one expects horizontal camera motion to be more common than vertical motion in hand-held and car-mounted videos. Thus the topology of the optical flow data set combines important patterns from both the range image and camera motion databases: the horizontal flow circle can be interpreted as the image of the range primary circle after applying common horizontal camera motions.

\FloatBarrier

\subsection{Remark}

Any method in data analysis must tolerate noisy data points, namely points $x \in X$ with low density. One may try to remove noisy points from the data set, but in general it is not easy to remove noise without also removing features. Because we use data set $X$ primarily to build density estimate $f \colon \R^n \to \R$, we believe that noise which is sparse enough to not significantly affect $f$ will not significantly affect our CW complex output. As evidence for this claim, consider the large samples $\cM$, $\cR$, and $\cF$ of optical image, range image, and optical flow patches. Both \cite{KleinBottle} and \cite{Range} remove noise from $\cM$ and $\cR$ by studying only core subsets of points with densities in the top 30\%. However, in this paper we work with random subsets of $\cM$ and $\cR$ containing both dense and noisy points. As shown in Figures~\ref{fig: OpticalImagePatches_1cells}, \ref{fig: RangeImagePatches_1cells}, and \ref{fig: OpticalFlowPatches_1cells}, we extract small CW complex models without first attempting to discard noisy data points.

\section{Conclusion}\label{S: Conclusion}

We have introduced a method for finding structure in high-dimensional Euclidean data sets. Following Morse theory, we use the nudged elastic band method to sample cells from the analogue of the Morse complex determined by the density function. We produce a nested family of CW complex models representing the dense regions of the data. We test the approach on social network, optical image, and range image data sets and find compact complexes revealing important nonlinear patterns. Moreover, we discover new topological structure, the horizontal flow circle, in a data set of optical flow patches.

\appendix
\appendixpage

\section{Additional data set information}\label{A: Additional data set information}

In this appendix we provide further information about each of the data sets. Our Java code and the optical image, range image, and optical flow data sets are available at the following webpage: \url{http://code.google.com/p/neb-tda}.

\subsection{Social network}

More information on the social network data set is available at the Adolescent Health Study webpage: \url{http://www.cpc.unc.edu/projects/addhealth}. To obtain a copy of the non-linkable Adolescent Health Network Structure files, please contact \url{addhealth@unc.edu}.

\subsection{Optical image patches}

The optical image database collected in \cite{VanHaterenVanDerSchaaf} contains 4,167 grayscale images from indoor and outdoor scenes, each $1020 \times 1532$ pixels (Figure~\ref{fig: opticalImages}). More details and the database itself are available at \url{http://bethgelab.org/datasets/vanhateren/}. From this database, Lee et al.\  in \cite{LeePedersenMumford} create a set $\cM$ of high-contrast, normalized, $3 \times 3$ patches through the following preprocessing steps.

\begin{enumerate}

\item Lee et al.\ select a large random sample of $3 \times 3$ patches with coordinates the logarithms of grayscale pixel values.
$$\begin{bmatrix}
x_1 & x_4 & x_7 \\
x_2 & x_5 & x_7 \\
x_3 & x_6 & x_9
\end{bmatrix}$$
Each patch is represented by a vector $x = (x_1, \dots, x_9)^T \in \R^9$. 

\item Lee et al.\ define a norm $\| \; \|_D$ measuring the contrast of a patch. Two coordinates $x_i$ and $x_j$ of $x$ are neighbors, denoted $i \sim j$, if the corresponding pixels in the $3 \times 3$ patch are adjacent. Let $\|x\|_D = \sqrt{\sum_{i \sim j} (x_i-x_j)^2}$. Lee et al.\ select the patches with contrast norm in the top 20\% of their sample.

\item Lee et al.\ normalize each patch by subtracting the average coordinate value and by dividing by the contrast norm. This maps the patches to a 7-dimensional ellipse.

\item Lee et al.\ change to the Discrete Cosine Transform (DCT) basis $\{e_1, \dots, e_8\}$ for $3 \times 3$ patches. The basis vectors are normalized to have contrast norm one, and so this maps the patches to a 7-dimensional sphere.

\end{enumerate}
Let $\cM$ be the resulting set of high-contrast, normalized, $3 \times 3$ optical patches. Three dense core subsets from $\cM$ are studied in \cite{KleinBottle} using persistent homology. In this paper we select a random data set $X \subset \cM$ of size 15,000.

\subsection{Range image patches}

The Brown range image database by Lee and Huang is a set of 197 range images from indoor and outdoor scenes, mostly $444 \times 1440$ pixels (Figure~\ref{fig: rangeImages}). The operational range for the Brown scanner is typically 2--200 meters, and the distance values for the pixels are stored in units of 0.008 meters. The database can be found at the following webpage: \url{http://www.dam.brown.edu/ptg/brid/index.html}.

From the Brown database we obtain a space of range image patches through the following steps, which are similar to the procedures in \cite{LeePedersenMumford, Range}.
\begin{enumerate}

\item We randomly select about $4 \cdot 10^5$ size $5 \times 5$ patches from the images in the database. Each patch is represented by a vector $x\in\R^{25}$ with logarithm values.

\item We compute the contrast norm $\| x \|_{D} =\sqrt{\sum_{i \sim j} (x_i-x_j)^2}$ of each patch and select the patches with contrast norm in the top 20\% of the entire sample.

\item We subtract from each patch the average of its coordinates and divide by the contrast norm.

\item We change to the DCT basis $\{e_1, \dots, e_{24}\}$ for $5 \times 5$ patches, normalized to have contrast norm one. This maps the patches to a 23-dimensional sphere.

\end{enumerate}
Let $\cR$ be the resulting set of high-contrast, normalized, $5 \times 5$ range patches. Our data set is a random subset $X \subset \cR$ of size 15,000.

\subsection{Optical flow patches}

The optical flow database by Roth and Black \cite{RothBlack} contains 800 optical flow fields, each $250\times 200$ pixels (Figure~\ref{fig: opticalFlowImages}). This ground-truth optical flow is generated by pairing range images from the Brown range image data base with camera motions. The camera motions are extracted from a database of 67 videos from hand-held or car-mounted video cameras, each approximately 100 frames long, using {\em boujou} software by 2d3 Ltd.\ available at \url{http://www.2d3.com}. The Roth and Black database is available at \url{http://www.gris.informatik.tu-darmstadt.de/~sroth/research/flow/downloads.html}.

From the Black and Roth database we create a space of optical flow patches, and our preprocessing is similar to that of \cite{LeePedersenMumford}.
\begin{enumerate}

\item We randomly choose $4\cdot10^5$ size $3\times3$ optical flow patches from the Roth and Black database. Each patch is a matrix of ordered pairs, where $u_i$ and $v_i$ are the horizontal and vertical components, respectively, of the flow vector at pixel $i$.
$$\begin{bmatrix}
(u_1,v_1) & (u_4,v_4) & (u_7,v_7) \\
(u_2,v_2) & (u_5,v_5) & (u_8,v_8) \\
(u_3,v_3) & (u_6,v_6) & (u_9,v_9)
\end{bmatrix}$$
We define $u = (u_1, \dots, u_9)^T$ and $v = (v_1,\dots, v_9)^T$ to be the vectors of horizontal and vertical flow components. We rearrange each patch to be a vector
$$x = \begin{pmatrix} u \\ v \end{pmatrix} \in \R^{18}.$$

\item We compute the contrast norm $\| x \|_D = \sqrt{\sum_{i \sim j} \| (u_i,v_i) - (u_j,v_j) \|^2}$ for all patches. We select the patches with contrast norm in the top 20\% of the entire sample.

\item We normalize the patches to have zero average flow. More explicitly, given a patch $x$, let $\overline{u}\in\R^9$ have each entry equal to $\frac{1}{9}\sum_{i=1}^{9}u_i$, the average of the horizontal components. Let $\overline{v}$ be defined similarly. We replace each patch $x$ with
$$\begin{pmatrix} u - \overline{u} \\ v - \overline{v} \end{pmatrix}.$$
We then divide each patch by its contrast norm.

\item We change to the DCT basis $\{e_1^u, \dots, e_8^u, e_1^v, \dots, e_8^v\}\subset\R^{18}$ for $3 \times 3$ optical flow patches, where
$$e_i^u =  \begin{pmatrix} e_i \\ \vec{0} \end{pmatrix} \quad \textrm{and} \quad e_i^v = \begin{pmatrix} \vec{0} \\ e_i \end{pmatrix}.$$
This maps the patches to a 16-dimensional sphere.

\end{enumerate}
Let $\cF$ be the resulting set of high-contrast, normalized, $3 \times 3$ optical flow patches. Our data set is a random subset $X \subset \cF$ of size 15,000.

\section{Initial bands}\label{A: Initial bands}

Let $p$ and $q$ be distinct 0-cells. Our method for generating the initial bands between $p$ and $q$ depends on whether data set $X$ is a general data set in $\R^n$ or whether $X$ is normalized to lie on a unit sphere $S^{n-1}\subset \R^n$.

For a general data set $X\subset\R^n$, we pick a uniformly random vector $y$ from the set of all unit vectors perpendicular to $p-q$, which is a sphere of dimension $n-2$. We also pick $r\in[0,d(p,q)]$ uniformly randomly. The resulting initial band is $N$ evenly distributed nodes along the circular arc (or straight line, with probability zero) between the points $p$, $(p+q+ry)/2$, and $q$.

If $X$ is normalized to lie on a unit sphere $S^{n-1}\subset \R^n$, we generate initial bands lying near $S^{n-1}$. Though $p$ and $q$ are near dense regions of $X$ they need not lie in $X$ nor in $S^{n-1}$. Let $\hat{p} = p/\|p\|$ and $\hat{q} = q/\|q\|$. We pick a uniformly random vector $y\neq \hat{p}, \hat{q}$ in $S^{n-1}$. The plane defined by $y$, $\hat{p}$, and $\hat{q}$ intersects $S^{n-1}$ in a circle. Let $\hat{p} = \hat{v}_1, \dots, \hat{v}_N = \hat{q}$ be the unique band that is evenly-spaced along this intersection circle, that starts at $\hat{p}$, that ends at $\hat{q}$, and that does not pass through $y$. We define our initial band to be $p = v_1, \dots, v_N = q$, where
$$v_i = \frac{(N-i) \|p\| + (i-1) \|q\|}{N-1}\ \hat{v}_i.$$

\section{Higher-dimensional cells}\label{A: Higher-dimensional cells}

One can imagine adapting the NEB method to search for higher-dimensional cells. Suppose for $k > 1$ that we have obtained the $(k-1)$-skeleton of a CW complex model. Now given an initial $k$-cell, we would like to move it towards a maximum density $k$-cell, and we describe a na\"ive approach. We model a $k$-cell by a graph $G = (V, E)$ that is the 1-skeleton of a regular CW complex homeomorphic to a $k$-cell. We choose an attaching map from the boundary of the $k$-cell to the $(k-1)$-skeleton of the CW complex model; this fixes the nodes in $V$ that lie in the boundary of the $k$-cell. Given a node $\alpha \in V$, let $v_\alpha\in \R^n$ denote its position and let $V_\alpha = \{ \beta \in V\ |\ \{\alpha, \beta\} \in E \}$ denote its set of incident vertices. We estimate the $k$-dimensional tangent space at $v_\alpha$ to be the span of the first $k$ components of a principal component analysis on $\{ v_\beta - v_\alpha\ |\ \beta \in V_\alpha \}$. If $\alpha$ is not a boundary node, then we define
$$F_\alpha = c \nabla f(v_\alpha)|_\perp\ +\ \sum_{\beta \in V_\alpha} (v_\beta - v_\alpha)$$
to be the force at vertex $\alpha$. As before, $\nabla f(v_\alpha)|_\perp$ is the component of $\nabla f(v_\alpha)$ perpendicular to the tangent space and is called the gradient force. Gradient constant $c$ adjusts the strength of the force. The term $\sum_{\beta \in V_\alpha} (v_\beta - v_\alpha)$ is the spring force. We numerically solve the system of first order differential equations $v_\alpha' = F_\alpha$.

We point out several weaknesses in the straightforward approach above.
\begin{itemize}
\item The spring force does not generalize the dimension $k=1$ case.
\item The gradient forces and spring forces do not both go to zero; instead, they balance against one another. This means that the cell does not converge exactly to the maximum density cell. A possible remedy is to project the spring force to the tangent space. One may then need to add an appropriate smoothing force to prevent kinks from forming in the cell. 
\item It may be preferable to have forces depend not only on the underlying 1-skeleton graph of the $k$-cell but also on its higher dimensional cells.
\item One may want an adaptive representation of a cell whose triangulation changes as the cell moves. Otherwise, the choice of an initial triangulation may affect the subsequent motion of the cell.
\item If the $(k-1)$-skeleton of the CW complex is complicated, then it is not clear which attaching maps to choose as the potential boundaries for $k$-cells.
\end{itemize}
We are interested in improved generalizations of NEB to higher-dimensional cells.

\begin{figure}[h]
\centerline{\includegraphics[width=.6\textwidth]{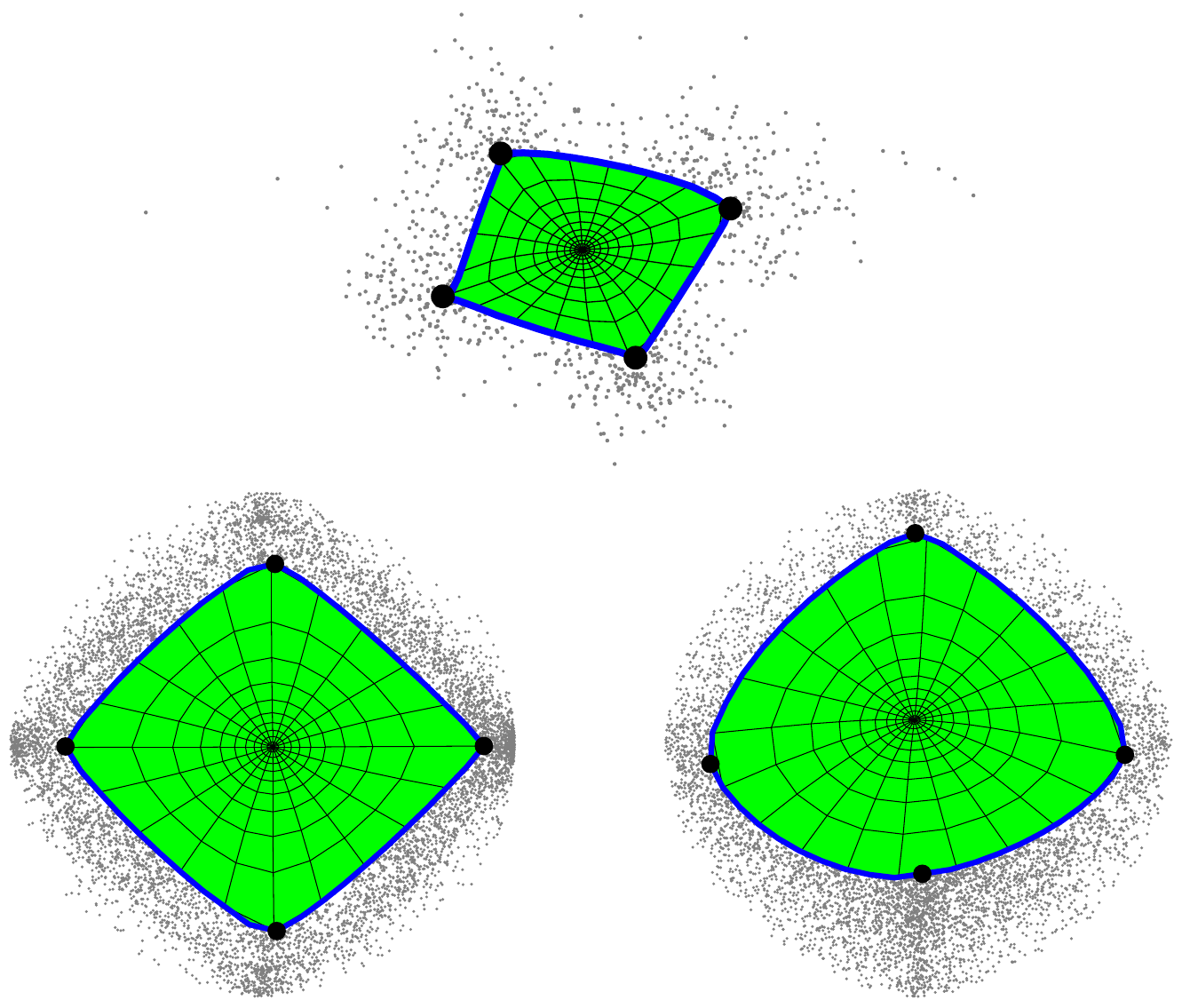}}
\caption{The 2-cells in the data sets. ({\it Top}) Social network data, projected to a plane using principal component analysis. ({\it Bottom left}) Range image patches, projected to the $e_1 e_5$ plane. ({\it Bottom right}) Optical flow patches, projected to the $e_1^u e_2^u$ plane.}
\label{fig: 2cells}
\end{figure}

Nevertheless, we use this approach to search for 2-cells in of our data sets. We use a particular CW complex homeomorphic to a 2-cell, containing 201 vertices, 400 edges, and 200 2-cells, whose 1-skeleton is a web-shaped graph with 20 nodes on each of 10 concentric rings (Figure~\ref{fig: 2cells}). We place the boundary nodes evenly-spaced on a chosen loop in the 1-skeleton; these nodes remain fixed. For the initial location of the 2-cell, we place the center node at the average of the boundary nodes, and we linearly interpolate between the boundary and center to place all remaining nodes. We estimate a convergent cell's density as $\min_{\alpha \in V} f(v_\alpha)$. See Figure~\ref{fig: 2cells} for the resulting 2-cells.

\begin{Acknowledgements}
We would like to thank Tim Harrington and Andrew Tausz for their help with the social network data and Guillermo Sapiro for his help with the optical flow data. The second author would like to thank the I.I.\ Rabi Science Scholars Program, Columbia University. This work is supported by Office of Naval Research Grant N00014-08-1-0931, Air Force Office of Scientific Research Grant FA9550-09-1-0143, Air Force Office of Scientific Research Grant FA9550-09-0-1-0531, and National Science Foundation Grant DMS-0905823.
\end{Acknowledgements}

\bibliographystyle{plain}

\bibliography{NudgedElasticBandInTopologicalDataAnalysis}

\end{document}